\documentclass[12pt]{article}

\begin{document}

\title{Limits of residually irreducible $p$-adic Galois representations}

\author{Chandrashekhar Khare}

\date{}

\maketitle

\newtheorem{theorem}{Theorem}
\newtheorem{lemma}{Lemma}
\newtheorem{prop}{Proposition}
\newtheorem{cor}{Corollary}
\newtheorem{conj}{Conjecture}
\newtheorem{guess}{Guess}
\newtheorem{remark}{Remark}
\newtheorem{example}{Example}
\newtheorem{conjecture}{Conjecture}
\newtheorem{definition}{Definition}
\newtheorem{quest}{Question}
\newtheorem{ack}{Acknowledgemets}
\newcommand{\rhobar}{\overline{\rho}}
\newcommand{\Sha}{{\rm III}}

\noindent{\bf Abstract:} In this note we produce examples
of converging sequences of Galois representations, and study
some of their properties.

\vspace{2mm}


\section{Introduction}

Consider a continuous representation $$\rho:G_L \to
GL_m(K)$$ of the absolute Galois group $G_L$ of 
a number field $L$,
with $K$ a finite extension of ${\bf Q}_p$, with $\cal O$ its
ring of integers, $| \ \ |$ its norm, and $k$ its residue field.
Then $\rho$ has an integral model taking values in
$GL_m({\cal O})$, and the semisimplification 
of its reduction modulo the maximal
ideal $\sf m$ of $\cal O$, denoted by $\rhobar$, 
is independent of the choice of
integral model. We assume that $\rhobar$ is absolutely irreducible and in fact will assume that all the $p$-adic 
representations considered in this paper are {\it residually absolutely 
irreducible}.

\begin{definition}\label{limits}
  An infinite sequence of (residually absolutely irreducible)
  continuous representations $\rho_i:G_{L}
  \rightarrow GL_m(K)$ tends to $\rho:G_{L} \to
  GL_m(K)$, if 
  $|{\rm tr}(\rho_i(g)) - {\rm tr}(\rho(g))| \rightarrow 0$ 
  uniformly for all $g \in G_{L}$. We also say that the $\rho_i$'s
  converge to $\rho$, or $\rho$ is their limit point.
\end{definition}
 

By Theorem 1 of [Ca], which we can apply because 
of our blanket assumption of residual absolute irreducibility,  
this is equivalent to saying that given any integer
$n$, for all $i>>0$, the reduction mod ${\sf m}^n$,
$\rho_{i,n}$, of (an integral model of) $\rho_i$
is isomorphic to the reduction mod ${\sf m}^n$, $\rho_n$,
of (an integral model of) $\rho$. 
Note that we are not assuming
that the $\rho_i$'s (or $\rho$) are {\it finitely ramified},
though we do know by the main theorem of [KhRa] that
the density of primes which ramify in a given $\rho_i$ is 0. 

In this note we study the limiting behavior of
the lifts produced in [R1] and 
completely characterise the limit points of these lifts (see
Theorem \ref{R-lifts} below). This suggests another approach to
certain special cases of the modularity lifting
theorems of Wiles, Taylor-Wiles et al. In the process 
we construct many sequences of converging $p$-adic
Galois representations (of fixed determinant and fixed ramification behaviour
at $p$). This raises many question 
that can be posed far more easily than answered.

Consider $\rhobar:G_{\bf Q} \rightarrow GL_2(k)$ that satisfies
the conditions of [R1], namely:

\begin{itemize}
      \item $\rhobar$ and ${\rm Ad}^0(\rhobar)$ are absolutely irreducible
            Galois representations, and the finite field $k$ 
            of characteristic $p$ is the minimal field of
            definition of $\rhobar$.
      \item The (prime to $p$) Artin conductor $N(\rhobar)$ of $\rhobar$ is 
            minimal amongst its twists. Denote by $S$ set of primes
            given by the union of the places where $\rhobar$ is
            ramified and $\{p,\infty\}$. 
      \item If $\rhobar$ is even then for the 
            decomposition group $G_p$ above $p$  we assume that $\rhobar 
            \mid_{G_p}$ is not twist equivalent to
            $\left(\begin{array}{cc} \chi &0 \\0 & 1\end{array}\right)$
            or twist equivalent to the indecomposable
            representation $\left(\begin{array}{cc} \chi^{p-2} & *\\
            0& 1\end{array}\right)$ where $\chi$ is the mod $p$
            cyclotomic character.
      \item If $\rhobar$ is odd we assume
             $\rhobar \mid_{G_p}$ is not
             twist equivalent to the trivial representation or the 
             indecomposable
             unramified representation given by
             $\left(\begin{array}{cc} 1 & * \\0 & 1\end{array}\right)$.
      \item $p\geq 7$ and the order of the projective image
            of $\rhobar$ is a multiple of $p$.
           
\end{itemize}

Let $Q=\{q_1,\cdots,q_n\}$ be a finite set of
primes such that $q_i \neq \pm 1$ mod $p$, unramified in $\rhobar$,
and the ratio of the eigenvalues of $\rhobar({\rm Frob}_{q_i})$
equal to $q_i^{\pm 1}$. We will call the primes as in $Q$ above 
{\it Ramakrishna primes for $\rhobar$} or {\it R-primes}
for short (suppressing the $\rhobar$ which is fixed). 
We consider the deformation ring $R_{S \cup
Q}^{Q-new}$ of [KR] (see Definition 1 of loc. cit.). To orient the 
reader we recall the definition of $R_{S \cup
Q}^{Q-new}$. For this we need:
\begin{definition}
  If $q$ is a prime, $G_{{\bf Q}_q}$ the absolute Galois group of ${\bf Q}_q$ and $R$ a complete
  Noetherian local ring with residue field $k$,  
  a continuous representation ${\rho}:G_{{\bf Q}_q} \rightarrow GL_2(R)$
  is said to be special if up to 
  conjugacy it is of the form $\left(\begin{array}{cc} \varepsilon\chi' & * \\0 & \chi'\end{array}\right)$ 
  for $\varepsilon$ the $p$-adic cyclotomic character,
  and $\chi': G_{{\bf Q}_q} \rightarrow R^*$ a continuous character.
  A continuous representation $\tilde{\rho}:G_{\bf Q} \rightarrow GL_2(R)$,
  is said to be special at a prime $q$ if $\tilde{\rho}|_{D_q}$, with $D_q$ a decomposition group at $q$,
  is special.
\end{definition}
Then $R_{S \cup Q}^{Q-new}$ is the universal ring that parametrises
deformations of $\rhobar$ that are minimally ramified 
at $S$ and such that at primes $q \in Q$ these deformations are 
{\it special}.  The ring
$R_{S \cup Q}^{Q-new}$ is a complete Noetherian local $W(k)$-algebra,
with $W(k)$ the Witt vectors of $k$.
The
deformation rings considered here are for the deformation problem with a
certain fixed (arithmetic) determinant character, 
and all the deformations of $\rhobar$
we consider will have this fixed determinant
character. 

\begin{definition}
  A finite set of R-primes $Q$ is said to be auxiliary
  if $R_{S \cup Q}^{Q-new} \simeq W(k)$.
\end{definition}

In [R1]
auxiliary sets  $Q$ of the above type 
were proven to exist. The representation
corresponding to 
$R_{S \cup Q}^{Q-new} \simeq W(k)$ 
is denoted by $\rho_{S \cup Q}^{Q-new}$. We will call these 
lifts {\it Ramakrishna lifts of $\rhobar$} or {\it R-lifts} 
for short (suppressing the $\rhobar$ which is fixed).

\begin{theorem}\label{R-lifts}
  A continuous representation $\rho:G_{\bf Q} \rightarrow
  GL_2(W(k))$ that is a deformation of $\rhobar$, is a limit point of distinct 
  R-lifts, if and only if $\rho$ is unramified
  outside $S$ and the set of all R-primes, and minimally ramified at 
  primes of $S$.
\end{theorem}

\noindent{\bf Remark:} Thus we have a complete description 
of the ``$p$-adic closure'' of {\it R-lifts}. Note that in particular
each $R$-lift is a limit point of other $R$-lifts. Note also
that any deformation $G_{\bf Q} \rightarrow GL_2(K)$ of $\rhobar$ that is
a limit point of $R$-lifts has a model that takes values in
$GL_2(W(k))$. The above theorem can be viewed in a sense as producing
an ``infinite fern'' structure (cf.,~[M]) in the set of all $R$-lifts of a given
$\rhobar$ as above. From the proof
of Theorem \ref{R-lifts} above, we in fact can deduce that
each $R$-lift gives rise
to infinitely many ``splines'' passing through it, where
a ``spline'' consists of a sequence of $R$-lifts converging to
it, and each element in a spline gives rise to its own
infinitely many splines. Missing from the picture are
the limit points of $R$-lifts which themselves are not $R$-lifts and 
which the theorem above characterises completely.

Here is the plan of the paper.
In Section 2 we prove Theorem \ref{R-lifts} which is a simple consequence of the methods of [R1] and [T1].
In Section 3 we prove a result about converging sequences of
representations arising from newforms, point
out a possible approach to the lifting theorems of Wiles et al that
is suggested by the work here. In Section 4 we
raise questions about rationality and motivic properties
of converging sequences of $p$-adic Galois representations. 

\section{Converging sequences of Galois representations}

We now prove Theorem \ref{R-lifts}.
It follows from
the methods of [R1] and [T1]. For the proof we need the following
lemma which follows from the methods of [R1] (see also Lemma 1.2 of [T1]) and Lemma 8 of [KR].

\begin{lemma}\label{borrowed} Let $\rho_n:G_{\bf Q} \rightarrow GL_2(W(k)/(p^n))$ 
be a lift of $\rhobar$ that is unramified
outside $S$ and the set of all R-primes, minimally ramified at 
primes of $S$, and special at all the primes outside $S$ at which it
is ramified.
Let $Q_n'$ be any finite set primes that includes the primes of ramification of $\rho_n$, 
such that $Q_n' \backslash S$ contains only $R$-primes and
such that $\rho_n|_{D_q}$ is special for $q \in Q_n' \backslash S$.
Then there exists a finite set of primes $Q_n$ that contains $Q_n'$,
such that $\rho_n|_{D_q}$ is special for $q \in Q_n \backslash S$, $Q_n \backslash S$ contains only $R$-primes
and $Q_n \backslash S$ is auxiliary. 
\end{lemma}

\noindent{\bf Proof:} We use [R1] and Lemma 8 of [KR] to construct an auxiliary set of primes $T_n$
such that $\rho_n|_{D_q}$ is special for $q \in T_n$. Then
as $Q_n' \backslash S$ contains only $R$-primes, it follows (using notation of [R1]) from Proposition 1.6 of [W]
that the kernel and cokernel of the map $$H^1(G_{S \cup T_n \cup Q_n'},{\rm Ad}^0(\rhobar))
\rightarrow \oplus _{v \in S \cup T_n \cup Q_n'}H^1(G_v,{\rm Ad}^0(\rhobar))/{\cal N}_v$$
have the same cardinality. Then using Proposition 10 of [R1],
or Lemma 1.2 of [T1], and Lemma 8 of [KR], we can augment the set $S \cup T_n \cup Q_n'$ to get a set $Q_n$ as in the statement of the lemma.
 
\vspace{3mm}

We are now ready to prove Theorem \ref{R-lifts}.
If $\rho:G_{\bf Q} \rightarrow
GL_2(W(k))$ is a limit point of  
$R$-lifts, then it is clear that $\rho$ is unramified
outside $S$ and the set of all $R$-primes, and minimally ramified at 
primes of $S$. We prove the converse. So let $\rho$ satisfy the
conditions of Theorem \ref{R-lifts}, and recall that we denote by
$\rho_n$ the reduction modulo $p^n$ of $\rho$.
It is easily checked that if $q$ is a $R$-prime any deformation of
$\rhobar|_{D_q}$ to a {\it ramified} $p$-adic representation 
is {\it special}: this follows
from the structure of tame inertia and the fact that $q^2 \neq 1$ mod $p$.
Further from the method of proof of Proposition 1 
of [KhRa], we easily deduce that the set of primes
$q$ for which $\rho|_{D_q}$ is special
is of density 0. Thus using Cebotarev and the assumptions on $\rhobar$ in
the introduction, we 
choose a finite set of primes $Q_n'$ such that
\begin{itemize}
\item $Q_n' \backslash S$ consists of $R$-primes and $\rho_n|_{D_q}$ is special for $q \in Q_n' \backslash S$,
\item $Q_n'$ contains all the ramified primes of $\rho_n$,
\item for some prime $q \in Q_n' \backslash S$, $\rho|_{D_q}$ is not special.
\end{itemize}
Using Lemma \ref{borrowed} 
we complete $Q'_n$ to a set $Q_n$ such that $Q_n \backslash S$ is auxiliary and
$\rho_n|_{D_q}$ is special for $q \in Q_n \backslash S$.
Then we claim $\rho_{S \cup Q_n}^{Q_n \backslash S-new} \equiv \rho$ mod $p^n$.
The claim is true as there is a unique representation $G_{\bf Q} \rightarrow
GL_2(W(k)/(p^n))$ (with the determinant that we have fixed)
that is unramified outside $S \cup Q_n$, minimal at $S$ and
special at primes of $Q_n \backslash S$ (as $R_{S \cup
Q_n}^{Q_n \backslash S-new} \simeq W(k)$). By construction the sets $Q_n$ contain at least one prime at
which $\rho$ is not special. Thus we see that we can pick a subsequence of
mutually distinct representations $\rho_i$ from
the $\rho_{S \cup Q_n}^{Q_n \backslash S-new}$'s such that $\rho_i \rightarrow
\rho$.

\vspace{3mm}

\noindent{\bf Remark:}
It is of vital importance that
$\rho$ is $GL_2(W(k))$-valued as otherwise we would not be 
able to invoke the disjointness results
that are used in the proof of Lemma \ref{borrowed} (Lemma 8 of [KR]). 

\vspace{3mm}

\noindent{\bf Remark:}
Theorem \ref{R-lifts}
can be applied in practise to give many examples of converging sequences
of $p$-adic representations:
for a non-CM elliptic curve $E_{/ {\bf Q}}$ for most primes $p$ 
the mod $p$ representation
satisfies the conditions given in the introduction, and the
corresponding $p$-adic representation is minimally ramified 
and $GL_2({\bf Z}_p)$ valued.

\vspace{3mm}

We end this section with a result that
refines the main result of [KhRa].

\begin{prop}\label{sparse}
  If $\rho_i:G_{L} \rightarrow GL_m(K)$ is a sequence of
  (residually absolutely irreducible) continuous representations that
  converges to $\rho$, then the set of primes where 
  any of the $\rho_i$'s is ramified (i.e.,
  $\cup {\rm Ram}(\rho_i)$ where ${\rm Ram}(\rho_i)$ is the set
  of primes at which $\rho_i$ is ramified) is
  of density zero.
\end{prop}

\noindent{\bf Proof:} Denote by $\rho_{i,n}$ (resp., $\rho_n$) 
the reduction mod ${\sf m}^n$
of an integral model of $\rho_i$ (resp., $\rho$). The proof 
consists in applying
Theorem 1 of [KhRa] twice: more precisely first its
statement, and then its proof. By an application of its statement we conclude
that the density
of $\cup_{i=1}^n{\rm Ram}(\rho_i)$ is 0 for any $n$. 
Now applying the proof of Theorem 1 of [KhRa], we 
define $c_{\rho,n}$ to be the upper density of the set
$S_{\rho,n}$
of primes $q$ of $L$ that
\begin{itemize}
\item lie above primes which split in $L/{\bf Q}$
\item are unramified in $\rho_1$ and $\neq p$,
\item $\rho_n|_{D_q}$ is unramified, but there exists a ``lift'' 
of $\rho_n|_{D_q}$, with $D_q$ the decomposition group at $q$, 
to a representation $\tilde{\rho}_q$ of $D_q$ 
to $GL_m(K)$ that is
ramified at $q$: by a lift we mean some conjugate of $\tilde{\rho}_q$
reduces mod ${\sf m}^n$ to $\rho_n|_{D_q}$.
\end{itemize}

We have from [KhRa] (see Proposition 1 of loc. cit. which was stated in greater
generality than needed there with the present application in mind): 

\begin{lemma}
  Given any $\varepsilon>0$, there is an integer $N_{\varepsilon}$ such that
  $c_{\rho,n} < \varepsilon$ for $n>N_{\varepsilon}$.
\end{lemma}

To prove Proposition \ref{sparse} 
it is enough to show that given any $\varepsilon>0$,
the upper density of the set $\cup {\rm Ram}(\rho_i)$ 
is $< \varepsilon$.
As $\cup_{i=1}^n{\rm Ram}(\rho_i)$ has density 0 for (the finite) 
$n$ that is the supremum
of the $i$'s such that $\rho_{i,N_{\varepsilon}}$ is not isomorphic
to $\rho_{N_{\varepsilon}}$, and $\rho_{N_{\varepsilon}}$ is finitely ramified,
it follows from the lemma above 
that the upper density
of $\cup {\rm Ram}(\rho_i)$  is $< \varepsilon$. Hence Proposition 
\ref{sparse}.

\vspace{3mm}

\noindent{\bf Remark:} One can ask for more refined 
information about the asymptotics of ramified 
primes in (limits of) residually absolutely irreducible $p$-adic Galois
representations. For instance
in Theorem 1 of [KhRa] one can ask (clued by Theorem 10 of [S1])
if the order of growth
of ramified primes can be proved to be bounded by $O(x^{1- {{1} \over
{2N}} +\epsilon})$, where $N$ 
is the $p$-adic analytic dimension of ${\rm im}(\rho)$, for any
$\epsilon>0$. Such quantitative refinements were asked 
for by Serre in an e-mail message to the author and 
seem hard and will require a new idea (that goes beyond [KhRa])
and a strong use of effective versions of the Cebotarev density theorem.

\section{Finite and infinite ramification}

Let $L$ be a number field and $K$ a finite extension of ${\bf Q}_p$
as before.

\begin{definition}
  We say that a residually absolutely irreducible continuous representation 
  $\rho:G_L \rightarrow GL_n(K)$ is motivic
  if $\rho$ arises as a subquotient of the $i^{\rm th}$ \'etale cohomology
  $H^i(X \times_L \overline{L},K)$ of a smooth projective 
  variety $X$ defined over a number
  field $L$.
\end{definition}

A {\it motivic} representation is {\it finitely ramified}.
In [R] examples of residually irreducible representations
$\rho:G_{\bf Q} \rightarrow GL_2(K)$ were constructed that were
infinitely ramified (see also the last section of [KR]).
Infinitely ramified $p$-adic representations cannot be
{\it motivic}. But they can arise as limits of $p$-adic
representations that are {\it motivic}. Fix an embedding 
$\overline{\bf Q} \rightarrow \overline{{\bf Q}_p}$.
Then as in [R] (and the last section of [KR]), 
there is a sequence of eigenforms 
$f_i \in S_2(\Gamma_0(N_i))$, for a sequence of squarefree integers
$N_i$ such that $N_i \rightarrow \infty$ and $(p,N_i)=1$, 
new of level $N_i$ such that the corresponding $p$-adic representations
$\rho_{f_i}:G_{\bf Q} \rightarrow GL_2({\bf Z}_p)$ have a
$p$-adic limit $\rho$, with $\rho$  {\it infinitely
ramified}. Such a $\rho$ is {\it non-motivic}, but is the limit
of {\it motivic} $p$-adic representations. Such limits of eigenforms
(in the works of
Serre and Katz for instance, cf., [Ka])
have been considered when varying weights
or varying the $p$-power level, while fixing the prime-to-$p$ part
of the level.

\begin{prop}\label{wiles}
  Let $f_i \in S_2(\Gamma_0(N_i))$
  be a sequence of eigenforms 
  with coefficients in a finite extension $K$
  of ${\bf Q}_p$ with $(N_i,p)=1$ and $p \geq 3$, that tend in the $p$-adic
  $q$-expansion topology to an element $f \in K[[q]]$, such the corresponding 
  residual representation $\rhobar$ satisfies the conditions
  in the introuction. The element
  $f$, that gives rise naturally to a Galois representation $\rho_f:G_{\bf Q}
  \rightarrow GL_2(K)$, is
  the $q$-expansion of a classical eigenform (of weight 2) if and only if $\rho_f$
  is finitely ramified.
\end{prop}

\noindent{\bf Proof:} The only if part is clear. The if part follows
from the methods of Wiles (see Chapter 3 of [W] and also [TW]) 
and their refinements: note that $\rho_f$ is finite, flat at $p$.

\vspace{3mm}

\noindent{\bf Remark:} 
Applying Theorem 1 when $\rhobar$ is odd
and finite flat at $p$, in which case
the $R$-lifts are modular by Theorem 1 of [K], we
can construct systematically many examples of sequences of eigenforms $f_i \rightarrow f$ ($f \in K[[q]]$),
with the levels of $f_i$ unbounded and such that $\rho_f$ is finitely
ramified ($f$ in fact then a classical eigenform as above). On the other hand as recalled above 
in [R] (see also last section of [KR]) we have examples of situations as above
with $\rho_f$ infinitely ramified.

\vspace{3mm}

It will be of interest to see if Proposition \ref{wiles}
could be proved in a more self-contained manner. The proof 
above does not use seriously the fact that one does know that $f$ arises
as a limit of the classical forms $f_i$. If such a proof could be
devised,
in conjunction with Theorem 1 above and Theorem 1 of [K] (which
is due to Ravi Ramakrishna) it would give in special cases 
a simpler approach using $R$-primes
to the modularity lifting theorems
of Wiles et al (see also [K]) that
directly works with the $p$-adic Galois representation 
that needs to be proved
modular, and if it could be implemented would avoid (in special cases albeit) the 
sophisticated deformation theoretic approach of [W]. 

We elaborate on this: Assume that $\rhobar$ is modular. 
In Theorem 1 we have characterised
the limit points of $R$-lifts. By Theorem 1 of [K] which proves that the
representation corresponding to $R_{S \cup Q}^{Q-new} \simeq W(k)$ is
modular as a consequence of the isomorphism 
$R_{S \cup Q}^{Q-new} \simeq {\bf T}_{S \cup Q}^{Q-new}$ (using notation of [K]),
we know that $R$-lifts are modular. Hence limits of $R$-lifts do
arise as limits of $p$-adic representations arising from
classical newforms. It only (!) remains to prove that a limit
of a converging sequence of 
$p$-adic representations arising from newforms (say of weight 2
and level prime to $p$ to avoid delicate considerations at $p$)
that is finitely ramified itself arises from a newform (i.e., prove 
Proposition \ref{wiles} without appealing directly to [W]).
Note that for a semistable elliptic curve $E$, for all large enough
primes $p$ (bigger than 3 for the methods here to directly work
unfortunately!), $T_p(E)$ is a limit point of $R$-lifts.

\vspace{3mm}

\noindent{\bf Note:} {\it In recent work we have indeed been able to
  give a self-contained approach to a result like Proposition
  \ref{wiles} above under some technical restrictions: see [K1].}

\section{Questions}

Proposition \ref{wiles} suggests that a representation that arises as a limit of motivic representations
(of ``bounded weights'': see Definition \ref{weight} below)
is finitely ramified if and only if it is motivic. We first recall one of the main conjectures in [FM]
in a form that is most pertinent for the considerations here.

\begin{conjecture}(Fontaine-Mazur)
  Consider a continuous 
  residually absolutely irreducible representation $\rho:G_L \rightarrow GL_m(K)$ that is
  potentially semistable at places above $p$. Then
  the following are equivalent:
\begin{enumerate}
\item $\rho$ is motivic
\item $\rho$ is finitely ramified.
\end{enumerate}
\end{conjecture}

From our earlier considerations it is natural to ask the following weaker question.

\begin{quest}\label{two} 
  Consider a continuous 
  residually absolutely irreducible representation $\rho:G_L \rightarrow GL_m(K)$ that is
  potentially semistable at places above $p$ that arises as the limit of motivic representations
  $\rho_i$. Then if $\rho$ is finitely ramified, is $\rho$ motivic? 
\end{quest}

\vspace{3mm}

It seems unlikely that the infinitely ramified representations produced in [R] are {\it algebraic}
(see definition below). This motivates the following considerations.

\begin{definition}\label{algebraic}
  A continuous (residually absolutely irreducible) 
  representation $\rho:G_L \to
  GL_m(K)$ is said to be algebraic 
  if there is a number field $F$ such that the
  characteristic polynomial of $\rho({\rm Frob}_q)$
  has coefficients in the ring of integers of $F$
  for all primes $q$ which are unramified
  in $\rho$. The minimal such field is the field of definition of $\rho$.
\end{definition}

As by the main theorem of [KhRa], the set of primes
at which $\rho$ ramifies is of density 0, the definition above is
a sensible one.

\begin{definition}\label{weight}
  A continuous (residually absolutely irreducible)
  algebraic representation $\rho:G_{L} \to
  GL_m(K)$ is said to be of weight $\leq t$ ($t \in {\bf Z}$) if for
  primes $q$ that are unramified in $\rho$, any root $\alpha$ of the
  characteristic polynomial of $\rho({\rm Frob}_q)$ satisfies
  $|\iota(\alpha)| \leq |k_q|^{{t-1} \over 2}$ for any embedding $\iota:\overline{\bf
  Q} \rightarrow {\bf C}$, with $k_q$ the residue field at $q$.
\end{definition}

\begin{quest}\label{rationality}
  If $\rho_i:G_L \rightarrow GL_m(K)$ 
  is an infinite sequence of (residually absolutely irreducible)
  distinct algebraic representations, all of 
  weight $\leq t$ for some fixed integer $t$, converging 
  to $\rho:G_L \rightarrow GL_m(K)$,
  and $K_i$ the field of definition of $\rho_i$,
  does $[K_i:{\bf Q}] \rightarrow \infty$ as $i \rightarrow \infty$?
\end{quest}

\noindent{\bf Remark:}
\begin{itemize}
\item  It is observed in [R] (this is a remark of Fred Diamond) that in the
  situation of Question \ref{rationality} only finitely many of the
  $\rho_i$'s can arise from elliptic curves: this is a consequence of
  the Mordell conjecture which gives that suitable twists of the
  classical modular curves $X(p^n)$ for $n>>0$ have finitely many
  $L$-valued points for a given number field $L$.
 \item If Question \ref{rationality} has a
  negative answer, using Proposition \ref{sparse}, 
  we deduce that for a set of primes $\{r\}$ of
  density one, the characteristic polynomials of 
  $\rho_i({\rm Frob}_r)$ are eventually constant. Hence we deduce that
  the characteristic polynomials of 
  $\rho({\rm Frob}_r)$ are defined and integral
  over a fixed number field $F$, i.e., $\rho$ is {\rm algebraic} (in
  the case when $\rho$ is infinitely 
  ramified this is linked to the questions below).
\end{itemize} 

\begin{quest}
  Let $\rho:G_L \rightarrow GL_m(K)$ be a continuous, 
  residually absolutely irreducible representation that is
  potentially semistable at places above $p$. Then
  are the following equivalent:
\begin{enumerate}
\item $\rho$ is motivic
\item $\rho$ is finitely ramified
\item $\rho$ is algebraic?
\end{enumerate}
\end{quest}

In the question above, the equivalence of 1 and 2 is the Fontaine-Mazur conjecture
recalled above: the possible equivalence of 3 to 1 and 2
is the main thrust of the question. One
might even ask the stronger question: If $\rho:G_L \rightarrow
GL_m(K)$, a continuous, residually absolutely irreducible
representation, is algebraic, then is $\rho$ {\it forced} to be both 
finitely ramified, {\it and} potentially semistable at places
above $p$? All the questions of this section have a positive 
answer when $m=1$. 

\section{Acknowledgements} I thank Gebhard B\"ockle, Christophe
Breuil,  Dipendra Prasad and Ravi
Ramakrishna for helpful correspondence,
Kirti Joshi and C.~S.~Rajan for useful conversations,
Haruzo Hida, David Rohrlich, J-P.~Serre and the anonymous referee 
for helpful remarks on the manuscript.

\section{References}

\noindent [Ca] Carayol, H., {\it Formes modulaires et repr\'esentations 
galoisiennes avec valeurs dans un anneau local complet}, 
in {\it $p$-adic monodromy 
and the Birch and Swinnerton-Dyer conjecture}, 213--237, 
Contemp. Math., 165, AMS, 1994. 



\vspace{3mm}

\noindent [FM] Fontaine, J.-M.,  Mazur, B., {\it Geometric Galois 
representations}, Elliptic curves, modular forms, and Fermat's last 
theorem, Internat. Press, Cambridge (1995), 41--78 





\vspace{3mm}

\noindent [K] Khare, C., {\it On isomorphisms between deformation rings and
Hecke rings}, preprint available at {\sf http://www.math.utah.edu/\~{ }shekhar/papers.html}.

\vspace{3mm}

\noindent [K1] Khare, C., {\it Modularity of $p$-adic Galois
representations via $p$-adic approximations}, in preparation.

\vspace{3mm}

\noindent [KhRa] Khare, C., Rajan, C.~S.~, {\it The density
of ramified primes in semisimple $p$-adic Galois representations},
International Mathematics Research Notices  no. 12 (2001), 601--607.

\vspace{3mm}

\noindent [KR] Khare, C., Ramakrishna, R., {\it
Finiteness of Selmer groups and deformation rings}, preprint available at
{\sf http://www.math.utah.edu/\~{ }shekhar/papers.html}.

\vspace{3mm}

\noindent [Ka] Katz, N., {\it Higher congruences between modular
forms}, Annals of Math. 101 (1975), 332--367.



\vspace{3mm}

\noindent [R] Ramakrishna, R., {\it Infinitely ramified representations},
Annals of Mathematics 151 (2000), 793--815.

\vspace{3mm}

\noindent [R1] Ramakrishna, R., {\it Deforming Galois representations and the
conjectures of Serre and Fontaine-Mazur}, to appear in Annals of Math.



\vspace{3mm}

\noindent [S1] Serre, J-P., {\it  Quelques applications du theoreme de
densite de Chebotarev}, Collected Works, Vol. 3, 563--641.



\vspace{3mm}

\noindent [T1] Taylor, R., {\it On icosahedral Artin
representations II}, preprint.

\vspace{3mm}

\noindent[TW] Taylor, R., Wiles, A., 
{\it Ring-theoretic properties of certain Hecke
algebras}, Ann. of Math. (2) 141 (1995), 553--572. 

\vspace{3mm}

\noindent [W] Wiles, A., {\it Modular elliptic curves and 
Fermat's last theorem}, Ann. of Math. 141 (1995), 443--551. 

\vspace{3mm}

\noindent {\it Addresses of the author}: Dept of Math, Univ of Utah, 155 S 1400 E,
Salt lake City, UT 84112. e-mail address: shekhar@math.utah.edu
 
\noindent School of Mathematics, 
TIFR, Homi Bhabha Road, Mumbai 400 005, INDIA. 
e-mail addresses: shekhar@math.tifr.res.in

\end{document}